\documentclass{article}
\usepackage{amsfonts,amsmath,amssymb}
\usepackage{graphics, epsfig}
\usepackage{color}
\usepackage{appendix}
\usepackage{ulem}
\usepackage[makeroom]{cancel}
\allowdisplaybreaks
\newtheorem{proposition}{Proposition}[section]

\newtheorem{lemma}{Lemma}[section]

\newtheorem{remark}{Remark}[section]

\numberwithin{equation}{section}
\numberwithin{theorem}{section}
\numberwithin{proposition}{section}
\numberwithin{lemma}{section}
\numberwithin{remark}{section}
\newcommand{\noi}{\noindent}

\newcommand{\dsty}{\displaystyle}
\newcommand{\txty}{\textstyle}

\newcommand{\al}{\alpha}

\newcommand{\be}{\beta}

\newcommand{\gm}{\gamma}
\newcommand{\dl}{\delta}

\newcommand{\varep}{\varepsilon}

\newcommand{\vp}{\varphi}
\newcommand{\sig}{\sigma}

\newcommand{\om}{\omega}

\newcommand{\z}{\zeta}

\newcommand{\df}[1]{\buildrel\mbox{\small def}\over{#1}}

\newcommand{\rr}{\mathbb{R}}

\newcommand{\LRA}{\Longrightarrow}

\newcommand{\bbox}{\vrule height.6em width.6em 
depth0em} 

\newcommand{\loc}{\operatorname{loc}}

\newcommand{\dist}{\operatorname{dist}}

\newcommand{\pl}{\partial}

\newcommand{\intl}{\int\limits}

\def\Xint#1{\mathchoice
    {\XXint\displaystyle\textstyle{#1}}%
    {\XXint\textstyle\scriptstyle{#1}}%
    {\XXint\scriptstyle\scriptscriptstyle{#1}}%
    {\XXint\scriptscriptstyle\scriptscriptstyle{#1}}%
    \!\int}
\def\XXint#1#2#3{\setbox0=\hbox{$#1{#2#3}{\int}$}
    \vcenter{\hbox{$#2#3$}}\kern-0.5\wd0}
\def\bint{\Xint-}
\def\dashint{\Xint{\raise4pt\hbox to7pt{\hrulefill}}}
\def\dashiint{\bint\kern-0.15cm\bint}

\newcommand{\ovl}[3]{\int_{#1}^{#2}\kern-#3pt\raise4pt\hbox to7pt{\hrulefill}\ }

\newcommand{\ovll}[3]{\intl_{#1}^{#2}\kern-#3pt\raise4pt\hbox to7pt{\hrulefill}\ }

\newcommand{\tvl}[2]{\iint_{#1}\kern-#2pt\raise4pt\hbox to7pt{\hrulefill}\ }

\newcommand{\bye}{\end{document}}

\input harnack_mono.mac
\begin{document}
\title{$1$-Dimensional Harnack Estimates}
\author
{Fatma Gamze D\"{u}zg\"{u}n\\
Hacettepe University\\
06800, Beytepe, Ankara, Turkey\\
email: {\tt gamzeduz@hacettepe.edu.tr}
\and
Ugo Gianazza\\
Dipartimento di Matematica ``F. Casorati", 
Universit\`a di Pavia\\ 
via Ferrata 1, 27100 Pavia, Italy\\
email: {\tt gianazza@imati.cnr.it}
\and
Vincenzo Vespri\\
Dipartimento di Matematica e Informatica ``U. Dini"\\
Universit\`a di Firenze\\  
viale Morgagni 67/A, 50134 Firenze, Italy\\
email: {\tt vespri@math.unifi.it}}
\date{}
\maketitle
\vskip.4truecm
\centerline{\it Dedicated to the memory of our friend Alfredo Lorenzi}
\vskip.4truecm
\vskip.4truecm
\begin{abstract}
{Let $u$ be a non-negative super-solution to a $1$-dimensional 
singular parabolic equation of $p$-Laplacian type 
($1<p<2$). If $u$ is 
bounded below on a time-segment $\{y\}\times(0,T]$ by a positive number $M$,
then it has a power-like decay of order $\frac p{2-p}$ with respect to the space variable $x$ in $\rr\times[T/2,T]$.  This fact, stated 
quantitatively in Proposition~\ref{Prop:1:1}, 
is a ``sidewise spreading of positivity'' of solutions
to such singular equations, and can be considered as a 
form of Harnack inequality. The proof of such an effect  is based on geometrical ideas.}  
\end{abstract}
\section{Introduction}
Let $E=(\al,\be)$ and define $E_{-\tau_o,T}=E\times(-\tau_o,T]$, for $\tau_o$, $T>0$.
Consider the non-linear diffusion equation
\begin{equation}\label{Eq:proto}
u_t-(|u_x|^{p-2}u_x)_x=0,\qquad 1<p<2.
\end{equation}
A function
\begin{equation}  \label{Eq:def:1}
u\in C_{\loc}\big(-\tau_o,T;L^2_{\loc}(E)\big)\cap L^p_{\loc}\big(-\tau_o,T; W^{1,p}_{%
\loc}(E)\big)
\end{equation}
is a local, weak super-solution to \eqref{Eq:proto}, if
for every compact set $K\subset E$ and every sub-interval $[t_1,t_2]\subset
(-\tau_o,T]$
\begin{equation}  \label{Eq:def:2}
\int_K u\vp dx\bigg|_{t_1}^{t_2}+\int_{t_1}^{t_2}\int_K \big[-u\vp_t
+|u_x|^{p-2}u_x\,\vp_x\big]dxdt\ge0
\end{equation}
for all non-negative test functions
\begin{equation*}
\vp\in W^{1,2}_{\loc}\big(-\tau_o,T;L^2(K)\big)\cap L^p_{\loc}\big(-\tau_o,T;W_o^{1,p}(K)%
\big).
\end{equation*}
This guarantees that all the integrals in \eqref{Eq:def:2} are convergent. These equations are termed singular since, for $1<p<2$, 
the modulus of ellipticity $|u_x|^{p-2}\to\infty$ 
as $|u_x|\to 0$.
\vskip.2truecm
\begin{remark}
{\normalfont Since we are working with \emph{local} solutions, we consider the domain $E_{-\tau_o,T}=E\times(-\tau_o,T]$, instead of dealing with the more natural $E_T=E\times(0,T]$,
in order to avoid problems with the initial conditions. The only role played by $\tau_o>0$ is precisely to get rid of any difficulty at $t=0$, and its precise value plays no role in the argument to follow.}
\end{remark}
\begin{proposition}\label{Prop:1:1} 
Let $u$ be a non-negative, local, weak super-solution to 
\eqref{Eq:proto} in $E_{-\tau_o,T}$, in the 
sense of \eqref{Eq:def:1}--\eqref{Eq:def:2}, satisfying 
\begin{equation}\label{Eq:0:1}
u(y,t)>M\ \ \ \forall t\in(0,\frac T2]
\end{equation}
for some $y\in E$, and for some $M>0$. 
Let $\bar\rho\df=\left(\frac{2^{2-p}T}{M^{2-p}}\right)^{\frac1p}$,
take $\rho\ge4\bar\rho$,
and assume that
\begin{equation*}
{B_{\rho}}(\bar x)\subset B_{4\rho}(y)\subset E,\qquad\text{where}\ \ \dist(\bar x,y)=2\rho.
\end{equation*} 
There exists $\bar\sig\in(0,1)$, that can be determined a priori, quantitatively  only 
in terms of the data, and independent of $M$ and $T$, such that 
\begin{equation}\label{Eq:0:5}
u(x,t)\ge\bar\sig M\left(\frac{\bar\rho}{\rho}\right)^{\frac p{2-p}}\quad\text{ for all }\> 
(x,t)\in B_{\rho}(\bar x)\times[\frac T4,\frac T2]
\end{equation} 
\end{proposition}
\subsection{Novelty and Significance}\label{SS:novelty}
The measure theoretical information on the ``positivity set'' in $\{y\}\times(0,\frac T2]$  
implies that such a positivity set 
actually ``expands'' sidewise in $\rr\times[\frac T4,\frac T2]$,
with a power-like decay of order $\frac p{2-p}$ with respect to the space variable $x$.  
Although considered a sort of natural fact, to our knowledge this result has never been proven before; it is the analogue of the power-like decay of order $\frac 1{p-2}$ with respect to the time variable $t$, known in the degenerate setting $p>2$ (see \cite{DBGV-acta}, \cite[Chapter~4, Section~4]{DBGV-mono}, \cite{kuusi}).
As the $t^{-\frac1{p-2}}$-decay is at the heart of the Harnack 
estimate for $p>2$, so Proposition~\ref{Prop:1:1} could be used to give a more streamlined proof of the Harnack inequality in the singular, super-critical range $\frac{2N}{N+1}<p<2$. This will be the object of future work, where we plan to address the general $N$-dimensional case. 

The proof is based on geometrical ideas, originally introduced in two different contexts: the energy estimates of \S~\ref{S:energy} and the decay of \S~\ref{S:decay} rely on a method introduced in \cite{LS-2009} in order to prove the H\"older continuity of solutions to an anisotropic elliptic equation, and further developed in \cite{DMV1,DMV2}; the change of variable used in the actual proof of Proposition~\ref{Prop:1:1} was used in \cite{DBGV-PAMS}.
\subsection{Further Generalization}
Consider partial differential equations of the form
\begin{equation}\label{Eq:1:1}
u_t-(\mathbf A(x,t,u,u_x))_x=0\quad\text{ weakly in } E_{-\tau_o,T},
\end{equation}
where the function $\mathbb A: E_{-\tau_o,T}\times\rr\times\rr\to\rr$ 
is only assumed to be measurable and subject to
the structure condition
\begin{equation}\label{Eq:1:2}
\left\{
\begin{array}{ll}
{\mathbf A(x,t,u,u_x)\,u_x\ge C_o|u_x|^p}\\
{|\mathbf A(x,t,u,u_x)|\le C_1|u_x|^{p-1}}
\end{array}\quad\text{ a.e. in }\> E_{-\tau_o,T},\right.
\end{equation}
where $1<p<2$, $C_o$ and $C_1$ are given positive constants. It is not hard to show that Proposition~\ref{Prop:1:1} holds also for weak super-solutions to \eqref{Eq:1:1}--\eqref{Eq:1:2}, since our proof is entirely based on the structural properties of \eqref{Eq:proto}, and the explicit dependence on $u_x$ plays no role. However, to keep the exposition simple, we have limited ourselves to the prototype case.
\section{Energy Estimates}\label{S:energy}
Let $u$ be a non-negative bounded, weak super-solution in $E_{-\tau_o,T}$, assume
\[
0\le\mu_-=\inf_{E_{-\tau_o,T}} u,
\]
and let $\om$ be a positive parameter. Without loss of generality we may assume that $0\in(\al,\be)$. For $\rho$ sufficiently small, so that $(-\rho,\rho)\subset(\al,\be)$, let 
\begin{align*}
&B_\rho=(-\rho,\rho),\qquad Q=B_\rho\times(0,T],\\
&B_\rho(y)=(y-\rho,y+\rho),\qquad Q(y)=B_\rho(y)\times(0,T],\\
&a\in(0,1),\ \ H\in(0,1]\quad\text{ parameters that will be fixed in the following},\\
&A=\{(x,t)\in Q(y):\ u(x,t)<\mu_-+(1-a)H\om\}.
\end{align*}
\begin{proposition}\label{Prop:2:1}
Let $u$ be a non-negative, local, weak super-solution to 
\eqref{Eq:proto} in $E_{-\tau_o,T}$, in the 
sense of \eqref{Eq:def:1}--\eqref{Eq:def:2}. 
There exists a positive constant $\gm=\gm(p)$, such that 
for every cylinder $Q(y)=B_\rho(y)\times(0,T]\subset E_{-\tau_o,T}$, 
and every piecewise smooth, cutoff 
function $\z$ vanishing on $\pl B_\rho(y)$, 
such that $0\le\z\le1$, and $\z_t\le0$,  
\begin{equation}\label{Eq:1:7}
\begin{aligned}
&\int_{B_\rho(y)\cap (A\times\{0\})}\left[\frac{(u-\mu_-+a\om H)^{2-p}}{2-p}-\frac{u}{(\om H)^{p-1}}\right]\z^p(x,0)\,dx\\
&+\iint_{Q(y)\cap A}\frac{|u_x|^p}{(u-\mu_-+a\om H)^{p}}\z^p\,dxdt
\le\gm\iint_{Q(y)\cap A}|\z_x|^p\,dxdt\\
&+\gm\iint_{Q(y)\cap A}(u-\mu_-+a\om H)^{2-p}\z^{p-1}|\z_t|\,dxdt.
\end{aligned}
\end{equation}
\end{proposition}
\noi{\it Proof} - Without loss of generality, we may assume $y=0$.
In the weak formulation of \eqref{Eq:proto}
take $\dsty\vp=G(u)\z^p$ as test function, with
\begin{equation*}
G(u)=\left[\frac1{(u-\mu_-+a\om H)^{p-1}}-\frac1{(\om H)^{p-1}}\right]_+,
\end{equation*}
and $\z=\z_1(x)\z_2(t)$, where $\z_1$ vanishes outside 
$B_{\rho}$ and satisfies
\begin{equation}\label{Eq:1:5} 
0\le\z_1\le1,\qquad \z_1=1\>\text{ in }\> B_{\frac \rho2},\qquad 
|\partial_x\z_1|\le\frac{\gm_1}{\rho},
\end{equation}
for an absolute constant $\gm_1$  independent of $\rho$, and $\z_2$ is monotone decreasing, and satisfies
\begin{equation}\label{Eq:1:6} 
0\le\z_2\le1,\qquad \z_2=1\>\text{ in }\> (0,\frac T2],\qquad\z_2=0\>\text{ for }\> t\ge T,\qquad |\partial_t\z_2|\le\frac{\gm_2}{T},
\end{equation}
for an absolute constant $\gm_2$  independent of $T$. It is easy to see that we have
\[
G'(u)=-\frac{p-1}{(u-\mu_-+a\om H)^{p}}\,\chi_{A}.
\]
Modulo a Steklov averaging process, we have
\begin{align*}
&\iint_Q u_t G(u)\z^p\,dxdt\\
&+\iint_Q\z^p G'(u)|u_x|^p\,dxdt+p\iint_Q G(u)\,|u_x|^{p-2}\z^{p-1}u_x\cdot \z_x dxdt\ge0,\\
&\\
&(p-1)\iint_{Q\cap A}\frac{|u_x|^p}{(u-\mu_-+a\om H)^{p}}\z^p\,dxdt\\
&\le p\iint_{Q\cap A}\z^{p-1}\frac{|u_x|^{p-1}}{(u-\mu_-+a\om H)^{p-1}}|\z_x|\,dxdt\\
&+\iint_{Q\cap A}\frac{u_t}{(u-\mu_-+a\om H)^{p-1}}\z^p\,dxdt-\iint_{Q\cap A}\frac{u_t}{(\om H)^{p-1}}\z^p\,dxdt,\\
&\\
&(p-1)\iint_{Q\cap A}\frac{|u_x|^p}{(u-\mu_-+a\om H)^{p}}\z^p\,dxdt\\
&\le p\iint_{Q\cap A}\z^{p-1}\frac{|u_x|^{p-1}}{(u-\mu_-+a\om H)^{p-1}}|\z_x|\,dxdt\\
&+\iint_{Q\cap A}\partial_t\left[\frac{(u-\mu_-+a\om H)^{2-p}}{2-p}-
\frac{u}{(\om H)^{p-1}}\right] \z^p\,dxdt,\\
&\\
&(p-1)\iint_{Q\cap A}\frac{|u_x|^p}{(u-\mu_-+a\om H)^{p}}\z^p\,dxdt\\
&\le p\iint_{Q\cap A}\z^{p-1}\frac{|u_x|^{p-1}}{(u-\mu_-+a\om H)^{p-1}}|\z_x|\,dxdt\\
&+\int_{B_\rho\cap(A\times\{T\})}\left[\frac{(u-\mu_-+a\om H)^{2-p}}{2-p}-
\frac{u}{(\om H)^{p-1}}\right] \z^p(x,T)\,dx\\
&-\int_{B_\rho\cap(A\times\{0\})}\left[\frac{(u-\mu_-+a\om H)^{2-p}}{2-p}-
\frac{u}{(\om H)^{p-1}}\right] \z^p(x,0)\,dx\\
&-p\iint_{Q\cap A}\left[\frac{(u-\mu_-+a\om H)^{2-p}}{2-p}-
\frac{u}{(\om H)^{p-1}}\right] \z^{p-1}\z_t\,dxdt.
\end{align*}
The second term on the right-hand side vanishes, as $\z(x,T)=0$. Therefore, an application of Young's inequality yields
\begin{align*}
&(p-1)\iint_{Q\cap A}\frac{|u_x|^p}{(u-\mu_-+a\om H)^{p}}\z^p\,dxdt\\
&+\int_{B_\rho\cap(A\times\{0\})}\left[\frac{(u-\mu_-+a\om H)^{2-p}}{2-p}-
\frac{u}{(\om H)^{p-1}}\right] \z^p(x,0)\,dx\\
&\le \frac{p-1}2\iint_{Q\cap A}\frac{|u_x|^p}{(u-\mu_-+a\om H)^{p}}\z^p\,dxdt+
\gm\iint_{Q\cap A}|\z_x|^p\,dxdt\\
&+p\iint_{Q\cap A}\frac{(u-\mu_-+a\om H)^{2-p}}{2-p}\z^{p-1}|\z_t|\,dxdt,
\end{align*}
where we have taken into account that $\z_t\le0$.
Therefore, we conclude
\begin{equation}
\begin{aligned}
&\int_{B_\rho\cap (A\times\{0\})}\left[\frac{(u-\mu_-+a\om H)^{2-p}}{2-p}-\frac{u}{(\om H)^{p-1}}\right]\z^p(x,0)\,dx\\
&+\frac{p-1}2\iint_{Q\cap A}\frac{|u_x|^p}{(u-\mu_-+a\om H)^{p}}\z^p\,dxdt
\le\gm\iint_{Q\cap A}|\z_x|^p\,dxdt\\
&+\gm\iint_{Q\cap A}(u-\mu_-+a\om H)^{2-p}\z^{p-1}|\z_t|\,dxdt,
\end{aligned}
\end{equation}
where the first term on the left-hand side is non--negative, since for $0<s<1-a<1$ the function $f(s)=\frac{(s+a)^{2-p}}t$ is monotone decreasing, and $f(1-a)=\frac1{1-a}>1$. \hfill\bbox
\begin{remark}
{\normalfont The constant $\gm$ deteriorates, as $p\to1$.}
\end{remark}
\begin{remark}
{\normalfont Even though in the next Section $H$ basically plays no role, we chose to state the previous Proposition with an explicit dependence also on $H$ for future applications.}
\end{remark}
\section{A Decay Lemma}\label{S:decay}
Without loss of generality, we may assume $\mu_-=0$. Let $M=\om$, $L\le\frac M2$, and suppose that 
\begin{equation}\label{Eq:1:9}
u(0,t)>M\ \ \ \forall t\in(0,\frac T2].
\end{equation}
Now, let $s_o$ be an integer to be chosen, and define
\begin{align*}
A_{s_o}&=\{t\in(0,\frac T2]:\ \exists\, x\in B_{\frac \rho2},\ u(x,t)<\frac L{2^{s_o}}\}\\
A(t)&=\{x\in B_{\frac \rho2}:\ u(x,t)<L(1-\frac 1{2^{s_o}})\},\qquad t\in(0,\frac T2]\\
A&=\left\{(x,t)\in B_{\frac \rho2}\times(0,\frac T2]:\ u(x,t)<L(1-\frac 1{2^{s_o}})\right\}.
\end{align*}
\begin{lemma}\label{Lm:3:1}
Let $u$ be a non-negative, local, weak super-solution to 
\eqref{Eq:proto} in $E_{-\tau_o,T}$, in the 
sense of \eqref{Eq:def:1}--\eqref{Eq:def:2}.
Let \eqref{Eq:1:9} hold and take 
\[
L\le\min\{\frac M2,\left(\frac T{\rho^p}\right)^{\frac1{2-p}}\}.
\] 
Then, for any $\nu\in(0,1)$, there exists a positive integer $s_o$ such that
\[
|\{t\in(0,\frac T2]:\ \exists\,x\in B_{\frac \rho2},\ u(x,t)\le\frac L{2^{s_o}}\}|\le\nu|(0,\frac T2]|.
\]
\end{lemma}
\noi{\it Proof} - Take $t\in A_{s_o}$: by definition,
there exists $\bar x\in B_{\frac \rho2}$ such that $u(\bar x,t)<L/{2^{s_o}}$. 
On the other hand, by assumption $\dsty u(0,t)>2L$, and therefore, $u(0,t)+(L/2^{s_o})>L$.
Hence
\[
\ln_+\frac{u(0,t)+\frac L{2^{s_o}}}{u(\bar x,t)+\frac L{2^{s_o}}}>(s_o-1)\ln 2,
\]
and we obtain
\begin{align*}
(s_o-1)\ln 2&\le\ln_+\left(\frac L{u(\bar x,t)+\frac L{2^{s_o}}}\right)-\ln_+\left(\frac L{u(0,t)+\frac L{2^{s_o}}}\right)\\
&=\int_0^{\bar x}\frac{\partial}{\partial x}\left(\ln_+\left(\frac L{u(\xi,t)+\frac L{2^{s_o}}}\right)\right)\,d\xi\\
&\le\int_{-\frac \rho2}^{\frac \rho2}\left|\frac{\partial}{\partial x}\left(\ln_+\left(\frac L{u(x,t)+\frac L{2^{s_o}}}\right)\right)\right|\,dx\\
&=\int_{B_{\frac\rho2}\cap A(t)}\left|\frac{\partial}{\partial x}\left(\ln_+\left(\frac L{u(x,t)+\frac L{2^{s_o}}}\right)\right)\right|\,dx.
\end{align*}
If we integrate with respect to time over the set $A_{s_o}$, we have
\begin{align*}
(s_o-1)|A_{s_o}|&\ln 2\le\int_0^{\frac T2}\int_{B_{\frac\rho2}\cap A(t)}\left|\frac{\partial}{\partial x}\left(\ln_+\left(\frac L{u(x,t)+\frac L{2^{s_o}}}\right)\right)\right|\,dx\\
&\le\left[\int_0^{\frac T2}\int_{B_{\frac\rho2}\cap A(t)}\left|\frac{\partial}{\partial x}\left(\ln_+\left(\frac L{u(x,t)+\frac L{2^{s_o}}}\right)\right)\right|^p\,dxdt\right]^{\frac1p}|Q|^{\frac{p-1}p}\\
&\le\left[\iint_{Q\cap A}\frac {|u_x|^p}{(u+\frac L{2^{s_o}})^p}\z^p\,dxdt\right]^{\frac1p}|Q|^{\frac{p-1}p}.
\end{align*}
Apply estimates \eqref{Eq:1:7} with $\dsty a=\frac1{2^{s_o}}$, $H\om=HM=L$.  The requirement $H\le1$ is satisfied, since $L\le\frac M2$. 
They yield
\begin{align*}
(s_o-1)|A_{s_o}|&\le\gm |Q|^{\frac{p-1}p}\left[\iint_{Q\cap A}|\z_x|^pdxdt\right]^{\frac1p}\\
&+\gm |Q|^{\frac{p-1}p}\left[\iint_{Q\cap A}(u+\frac L{2^{s_o}})^{2-p}|\z_t|dxdt\right]^{\frac1p}.
\end{align*}
With these choices, we have
\begin{align*}
(s_o-1)|A_{s_o}|&\le\gm |Q|^{\frac{p-1}p}\frac1\rho|Q|^{\frac{1}p}+\gm|Q|^{\frac{p-1}p}\left(\frac{L^{2-p}}T\right)^{\frac1p}|Q|^{\frac{1}p}\\
&\le\gm\left[\frac1\rho+\left(\frac{L^{2-p}}T\right)^{\frac1p}\right]|Q|.
\end{align*}
If we require $\dsty L\le\left(\frac T{\rho^p}\right)^{\frac1{2-p}}$, and we substitute it back in the previous estimate, we have
\[
(s_o-1)|A_{s_o}\le\gm_1|(0,\frac T2]|.
\]
Therefore, if we want that $\dsty |A_{s_o}|\le\nu|(0,\frac T2]|$, it is enough to require that $\dsty s_o=\frac{\gm_1}\nu+1$. \hfill\bbox
\vskip.2truecm
\noi The previous result can also be rewritten as
\begin{lemma}\label{Lm:3:2}
Let $u$ be a non-negative, local, weak super-solution to 
\eqref{Eq:proto} in $E_{-\tau_o,T}$, in the 
sense of \eqref{Eq:def:1}--\eqref{Eq:def:2}.
Let \eqref{Eq:1:9} hold. For any $\nu\in(0,1)$, there exists a positive integer $s_o$ such that
\[
|\{t\in(0,\frac T2]:\ \exists\,x\in B_{\frac \rho2},\ u(x,t)\le\left(\frac T{\rho^p}\right)^{\frac1{2-p}}\frac1{2^{s_o}}\}|\le\nu|(0,\frac T2]|,
\]
provided $\rho>0$ is so large that $\dsty \left(\frac T{\rho^p}\right)^{\frac1{2-p}}\le\frac M2$.
\end{lemma}
Now let $\bar \rho$ be such that 
\begin{equation}\label{Eq:rhobar}
\left(\frac T{\bar \rho^p}\right)^{\frac1{2-p}}=\frac M2\qquad\Rightarrow\qquad\bar \rho=\left(\frac{2^{2-p}T}{M^{2-p}}\right)^{\frac1p},
\end{equation}
and assume that $B_{\bar \rho}\subset(\al,\be)$.
Then Lemmas~\ref{Lm:3:1}--\ref{Lm:3:2} can be rephrased 
\begin{lemma}\label{Lm:3:3}
Let $u$ be a non-negative, local, weak super-solution to 
\eqref{Eq:proto} in $E_{-\tau_o,T}$, in the 
sense of \eqref{Eq:def:1}--\eqref{Eq:def:2}.
Let \eqref{Eq:1:9} hold. For any $\nu\in(0,1)$, there exists a positive integer $s_o$ such that for any $\rho>\bar \rho$
\[
|\{t\in(0,\frac T2]:\ \exists\,x\in B_{\frac \rho2},\ u(x,t)\le\frac M{2^{s_o+1}}\left(\frac{\bar \rho}{\rho}\right)^{\frac p{2-p}}\}|\le\nu|(0,\frac T2]|,
\]
provided that $B_\rho\subset(\al,\be)$.
\end{lemma}
\begin{remark}
{\normalfont The previous corollary gives us the power-like decay, required in Proposition~\ref{Prop:1:1}.}
\end{remark}
Let us now set $A^c_{s_o}\df=(0,\frac T2]\backslash A_{s_o}$. Then, if \eqref{Eq:1:9} holds, we conclude that for any $t\in A^c_{s_o}$ and for any $x\in B_{\frac \rho2}$ with $\rho>\bar \rho$
\begin{equation}\label{Eq:1:10}
u(x,t)\ge\frac{M}{2^{s_o+1}}\left(\frac{\bar \rho}\rho\right)^{\frac p{2-p}}.
\end{equation}
Let $c\ge4$ denote a positive parameter, choose $\bar x\in(\al,\be)$ such that $|\bar x|=2c\bar \rho$, and consider $B_{c\bar \rho}(\bar x)$. Then, by \eqref{Eq:1:10}
\begin{equation}\label{Eq:1:11}
\forall\,x\in B_{c\frac{\bar \rho}2}(\bar x),\ \ \forall\,t\in A^c_{s_o}\ \ u(x,t)\ge\frac{M}{2^{s_o+1}}\left(\frac2{5c}\right)^{\frac p{2-p}},
\end{equation}
provided \eqref{Eq:1:9} holds, and $B(\bar x)_{c\bar \rho}\subset(\al,\be)$.
\section{A DeGiorgi-Type Lemma}
Assume that some information is available on 
the ``initial data'' relative to the cylinder 
$B_{2\rho}(y)\times(s,s+\theta\rho^p]$, say for example 
\begin{equation}\label{Eq:2:6}
u(x,s)\ge M\quad\text{ for a.e. }\> x\in B_{2\rho}(y)
\end{equation}
for some $M>0$. Then, the following Proposition is proved in \cite[Chapter~3, Lemma~4.1]{DBGV-mono}.
\begin{lemma}\label{Lm:2:2} 
Let $u$ be a non-negative, local, weak super-solution to 
\eqref{Eq:proto}, and $M$ be a positive number such that 
\eqref{Eq:2:6} holds. 
Then 
\begin{equation*}
u\ge{\txty\frac12} M\quad\text{ a.e. in }\> 
B_\rho(y)\times(s,s+\theta(2\rho)^p],
\end{equation*}
where 
\begin{equation}\label{Eq:2:7} 
\theta=\dl\,{M^{2-p}}
\end{equation}
for a constant $\dl\in(0,1)$ depending only upon $p$, and independent of $M$ and $\rho$.
\end{lemma}
\begin{remark}\label{Rmk:2:2} 
{\normalfont Lemma~\ref{Lm:2:2} is based on the 
energy estimates and Proposition~3.1 of 
\cite{dibe-sv}, Chapter~I which continue to hold in 
a stable manner for $p\to1$. These results are therefore 
valid for all $p\ge1$, including a seamless transition 
from the singular range $p<2$ to the degenerate range $p>2$.
}
\end{remark}
\section{Proof of Proposition~\ref{Prop:1:1}}
Fix $y\in E$, define $\bar \rho$ as in \eqref{Eq:rhobar}, and choose a positive parameter $C\ge 4$, such that the cylindrical domain
\begin{equation}\label{Eq:3:1}
B_{2^{\frac{2-p}p}C\bar\rho}(y)\times\big(0,\frac T2\big]\subset E_{-\tau_o,T}.
\end{equation}
This is an assumption both on the size of the reference ball $B_{2^{\frac{2-p}p} C\bar\rho}(y)$ and on $T$: nevertheless, we can always assume it without loss of generality. Indeed, if \eqref{Eq:3:1} were not satisfied, we would decompose the interval $(0,\frac T2]$ in smaller subintervals, each of width $\tau$, such that \eqref{Eq:3:1} is satisfied working with $\bar\rho$ replaced by 
$$\rho^*=\left(\frac{2^{2-p}\tau}{M^{2-p}}\right)^{\frac1p}.$$ 
The only role of $C$ is in determining a sufficiently large reference domain 
$$B_{2^{\frac{2-p}p} C\bar\rho}(y)\subset E,$$
which contains the smaller ball we will actually work with, 
and will play no other role; in particular the structural constants will not depend on $C$.
 
\noi Now, introduce the change of variables and the new unknown function
\begin{equation}\label{Eq:3:2}
z=2^{\frac{2-p}p}\frac{x-y}{\bar\rho},\qquad 
-e^{-\tau}=\frac{t-\frac T2}{\frac T2},
\qquad v(z,\tau)=\frac1M u(x,t)e^{\frac{\tau}{2-p}}.
\end{equation}
This maps the cylinder in (\ref{Eq:3:1}) into $B_C\times(0,\infty)$ 
and transforms (\ref{Eq:proto}) into 
\begin{equation}\label{Eq:3:3}
v_\tau-(|v_z|^{p-2}v_z)_z=\frac1{2-p}v 
\quad\text{ weakly in }\> B_C\times(0,\infty).
\end{equation}
The assumption 
(\ref{Eq:0:1}) of Proposition~\ref{Prop:1:1} becomes 
\begin{equation}\label{Eq:3:5}
v(0,\tau)\ge e^{\frac{\tau}{2-p}}\quad\text{ for all }\> \tau\in(0,+\infty).
\end{equation}
Let $\tau_o>0$ to be chosen and set 
\begin{equation*}
k=e^{\frac{\tau_o}{2-p}}.
\end{equation*}
With this symbolism, \eqref{Eq:3:5} implies
\begin{equation}\label{Eq:3:7}
v(0,\tau)\ge k\quad\text{ for all }\> \tau\in(\tau_o,+\infty).
\end{equation}
Now consider the segment
\begin{equation*}
I\df=\{0\}\times\big(\tau_o,\tau_o+k^{2-p}).
\end{equation*}
Let $\nu=\frac14$ and $s_o$ be the corresponding quantity introduced in Lemma~\ref{Lm:3:1}.
We can then apply Lemmas~\ref{Lm:3:1}--\ref{Lm:3:3} with $T=k^{2-p}$, $M$ substituted by $k$,
\[
A_{s_o}=\{\tau\in(\tau_o,\tau_o+\frac12]:\ \exists\, z\in B_{\frac \rho2},\ v(z,\tau)<\frac k{2^{s_o+1}}\}\quad\text{ for }\ \ \rho>\tilde\rho,
\]
with $\tilde\rho\df=2^{\frac{2-p}p}$. Therefore, if $c\ge4$ denotes a positive parameter, we choose $\bar z\in B_C$ such that $|\bar z|=2c\tilde\rho$, and consider $B_{c\tilde\rho}(\bar z)$, by \eqref{Eq:1:10}
\begin{equation}\label{Eq:4:11}
\forall\,z\in B_{c\frac{\tilde\rho}2}(\bar z),\ \ \forall\,\tau\in A^c_{s_o}\ \ v(z,\tau)\ge\frac{k}{2^{s_o+1}}\left(\frac2{5c}\right)^{\frac p{2-p}},
\end{equation}
provided $B(\bar z)_{c\tilde\rho}\subset B_C$.
Summarising, there 
exists at least a time level $\tau_1$ in the range 
\begin{equation}\label{Eq:3:13}
\tau_o<\tau_1<\tau_o+\frac12 k^{2-p}
\end{equation}
such that 
\begin{equation*}
\forall\,z\in B_{c\frac{\tilde\rho}2}(\bar z),\ \ v(z,\tau_1)\ge\sig_o e^{\frac{\tau_o}{2-p}}
\quad\text{ where }\quad \sig_o= \frac{1}{2^{s_o+1}}\left(\frac2{5c}\right)^{\frac p{2-p}}.
\end{equation*}
\begin{remark}\label{Rmk:3:2} {\normalfont Notice that $\sig_o$ is determined only in terms 
of the data and is independent of the parameter $\tau_o$, which is still to be chosen.
}
\end{remark}
\subsection{Returning to the Original Coordinates}\label{S:3:3}
In terms of the original coordinates and the 
original function $u(x,t)$, this implies 
\begin{equation*}
u(\cdot,t_1)\ge\sig_o M e^{-\frac{\tau_1-\tau_o}{2-p}}
\,\df{=}\, M_o
\qquad\text{ in }\> B_{c\frac{\bar\rho}2}(\bar x)
\end{equation*}
where the time $t_1$ corresponding to $\tau_1$ is 
computed from (\ref{Eq:3:2}) and (\ref{Eq:3:13}), and
$\dsty\dist(\bar x,y)=2c\bar\rho$.  
Now, apply Lemma~\ref{Lm:2:2} with $M$ replaced 
by $M_o$ over the cylinder $B_{c\frac{\bar\rho}2}(\bar x)\times
\big(t_1, t_1+\theta(c\bar\rho)^p\big]$.
By choosing  
\begin{equation*}
\theta = \dl M_o^{2-p}\qquad\text{ where }\quad 
\dl=\dl(\text{data}), 
\end{equation*}
the assumption (\ref{Eq:2:7}) is satisfied, and 
Lemma~\ref{Lm:2:2} yields
\begin{equation}\label{Eq:3:16}
\begin{aligned}
 u(\cdot, t)&\ge{\txty\frac12} M_o
={\txty\frac12}\sig_o M e^{-\frac{\tau_1-\tau_o}{2-p}}\\
&\ge\frac{1}{2^{s_o+2}}\left(\frac2{5c}\right)^{\frac p{2-p}} e^{-\frac2{2-p}e^{\tau_o}} M
\end{aligned}
\qquad
\text{ in }\quad B_{c\bar\rho}(\bar x)
\end{equation}
for all times
\begin{equation}\label{Eq:3:17}
t_1\le t\le t_1+ \dl \frac{1}{2^{s_o(2-p)}}\left(\frac25\right)^p e^{-(\tau_1-\tau_o)}\frac T2. 
\end{equation}
Notice that \eqref{Eq:3:16}
can be rewritten as
\begin{equation}\label{Eq:3:16bis}
u(\cdot, t)\ge\bar\sig\left(\frac{\bar\rho}\rho\right)^{\frac p{2-p}} M
\text{ in }\quad B_{\rho}(\bar x),
\end{equation}
with
\begin{equation}\label{Eq:3:16ter}
\bar\sig\df=\frac{1}{2^{s_o+2}}\left(\frac25\right)^{\frac p{2-p}}e^{-\frac2{2-p}e^{\tau_o}}
\end{equation}
If the right hand side of \eqref{Eq:3:17} equals $\frac T2$, then (\ref{Eq:3:16}) holds for 
all times in 
\begin{equation}\label{Eq:3:18} 
\big(\frac T2-\varep M^{2-p}(c\bar\rho)^p\,,\, \frac T2\big]\qquad\text{ where }\quad 
\varep=\dl \sig_o^{2-p} e^{-e^{\tau_o}};
\end{equation} 
taking into account the expression for $\bar\rho$ and $\sig_o$, we conclude 
that (\ref{Eq:3:16}) holds for 
all times in 
\begin{equation}\label{Eq:3:18bis} 
\big(\frac T2-e^{-e^{\tau_o}}\frac{\dl}{2^{s_o(2-p)}}\left(\frac25\right)^p\frac T2\,,\, \frac T2\big].
\end{equation} 

Thus, the conclusion of Proposition~\ref{Prop:1:1} holds, 
provided the upper time level in (\ref{Eq:3:17}) equals $\frac T2$.  
The transformed $\tau_o$ level is still undetermined, and it 
will be so chosen as to verify such a requirement. 
Precisely, taking into account (\ref{Eq:3:2})
\begin{equation*} 
\frac T2 e^{-\tau_1}=-(t_1-\frac T2)
=\dl \frac{1}{2^{s_o(2-p)}}\left(\frac25\right)^p e^{-(\tau_1-\tau_o)}\frac T2\,\,\LRA\,\,
e^{\tau_o}=\left(\frac52\right)^p\frac{2^{s_o(2-p)}}{\dl}.
\end{equation*}
This determines quantitatively $\tau_o=\tau_o(\text{data})$, and  
inserting such a $\tau_o$ on the right-hand side 
of \eqref{Eq:3:16ter} and \eqref{Eq:3:18bis}, yields a bound below that depends only on the data;
\eqref{Eq:3:16ter} and \eqref{Eq:3:18bis} have been obtained relying on the bound below for $u$ along the segment $\{y\}\times(0,\frac T2]$. However, the same argument on the bound along the shorter segment $\{y\}\times(0,s]$ for any $\frac T4\le s<\frac T2$ yields the same result with $\frac T2$ substituted by $s$: the proof of Proposition~\ref{Prop:1:1} is then completed.\hfill\bbox
\subsection{A Remark about the Limit as $p\to2$}\label{S:4}
The change of variables (\ref{Eq:3:2}) and the subsequent 
arguments,  yield constants that deteriorate as $p\to2$. 
This is no surprise, as the decay of solutions to linear parabolic equations 
is not power-like, but rather exponential-like, as in the fundamental solution of the heat equation.

Nevertheless, our estimates can be stabilised, in order to recover the correct exponential decay in the 
$p=2$ limit. However, this would require a careful tracing of all the functional dependencies in our estimates, and we postpone it to a future work.


\begin{thebibliography}{99}
\bibitem{dibe-sv} E. DiBenedetto, {\it Degenerate Parabolic 
Equations}, Universitext, Springer--Verlag, New York, 1993.
\bibitem{DBGV-acta} E. DiBenedetto, U. Gianazza and V. Vespri, 
Harnack Estimates for Quasi-Linear Degenerate Parabolic 
Differential Equation, {\it Acta Mathematica, 200 (2008), 181--209.}
\bibitem{DBGV-mono} E. DiBenedetto, U. Gianazza and V. Vespri, 
{\it Harnack's Inequality for Degenerate and Singular Parabolic 
Equations}, Springer Monographs in Mathematics, Springer-Verlag, 
New York, 2012.
\bibitem{DBGV-PAMS}E. DiBenedetto, U. Gianazza and V. Vespri, A New Approach to the Expansion of Positivity Set of Non-negative Solutions to Certain Singular Parabolic Partial Differential Equations, {\it Proc. Amer. Math. Soc. 138 (2010), 3521--3529.}
\bibitem{DMV1} F.G. D\"uzg\"un, P. Marcellini and V. Vespri, An alternative approach to the Hoelder continuity of solutions to some elliptic equations, {\it NonLinear Anal. (2014) 133--141.}
\bibitem{DMV2} F.G. D\"uzg\"un, P. Marcellini and V. Vespri, Space expansion for a solution of an anisotropic $p$-Laplacian equation by using a parabolic approach, {\it Riv. Mat. Univ. Parma, 5, (2014) 93--111.}
\bibitem{kuusi} T. Kuusi, The weak Harnack estimate for weak 
supersolutions to nonlinear 
degenerate parabolic equations, {\it Ann. Scuola Norm. Sup. Pisa 
Cl. Sci. (5),  {\bf7}(4), (2008), 673--716.} 
\bibitem{LS-2009} V. Liskevich and I.I. Skrypnik, H\"older continuity of solutions to an anisotropic elliptic equation, {\it Nonlinear Anal. 71 (2009), 1699--1708.}
\end{thebibliography}
\bye